\documentclass{amsart}
\usepackage{amsmath,amssymb}
\usepackage{yhmath}

\newcommand{\bdis}{\begin{displaymath}}
\newcommand{\edis}{\end{displaymath}}
\newcommand{\be}{\begin{equation}}
\newcommand{\ee}{\end{equation}}
\newcommand{\mbb}{\mathbb}
\newcommand{\mcal}{\mathcal}

\newcommand{\vp}{\varphi}

\newcommand{\zf}{\zeta\left(\frac{1}{2}+it\right)}
\newcommand{\zfvp}{\zeta\left(\frac{1}{2}+i\vp_1(t)\right)}

\newcommand{\onT}{\overset{0}{T}}
\newcommand{\ojT}{\overset{1}{T}}

\DeclareMathOperator{\im}{Im}

\theoremstyle{definition}

\theoremstyle{remark}
\newtheorem{remark}[]{Remark}

\newtheorem*{mydef11}{{\bf Theorem 1}}

\newtheorem*{mydef21}{{\bf Definition 1}}

\newtheorem*{mydef41}{{\bf Corollary 1}}

\newtheorem*{mydef51}{{\bf Lemma 1}}

\newtheorem*{mydef81}{{\bf Property 1}}

\newtheorem*{mydef91}{{\bf Formula1}}

\newtheorem*{mydef92}{{\bf Formula2}}

\newtheorem*{mydef93}{{\bf Formula3}}

\newtheorem*{mydef94}{{\bf Formula4}}

\numberwithin{equation}{section}

\begin{document}

\title[Jacob's ladders and vector operator\dots]{Jacob's ladders and vector operator producing new generations of $L_2$-orthogonal systems connected with the Riemann's $\zf$ function}{}

\author{Jan Moser}

\address{Department of Mathematical Analysis and Numerical Mathematics, Comenius University, Mlynska Dolina M105, 842 48 Bratislava, SLOVAKIA}

\email{jan.mozer@fmph.uniba.sk}


\begin{abstract}
In this paper we introduce a generating vector-operator acting on the class of functions $L_2([a,a+2l])$. This operator produces (for arbitrarily fixed $[a,a+2l]$) infinite number of new generation $L_2$-systems. Every element of the mentioned systems depends on Riemann's zeta-function and on Jacob's ladder.  
\end{abstract}
\maketitle

\section{Introduction} 

\subsection{}   

In this paper we introduce vector operator $\hat{G}$ defined on the class of all $L_2$-orthogonal systems 
\bdis 
\{f_n(t)\}_{n=0}^\infty,\ t\in [a,a+2l],\ \forall a\in\mbb{R},\ \forall l\in\mbb{R}^+ 
\edis 
that for a fixed class $L_2([a,a+2l])$ associates following new classes: 
\be \label{1.1} 
\{f_n(t)\}_{n=0}^\infty\xrightarrow{\hat{G}}\{f_n^{p_1}(t)\}_{n=0}^\infty,\ p_1=1,\dots,k, 
\ee  
\be \label{1.2} 
\{f_n^{p_1}(t)\}_{n=0}^\infty\xrightarrow{\hat{G}}\{f_n^{p_1,p_2}(t)\}_{n=0}^\infty,\ p_1,p_2=1,\dots,k, 
\ee  
and so on up to 
\be \label{1.3} 
\{f_n^{p_1,\dots,p_{s-1}}(t)\}_{n=0}^\infty\xrightarrow{\hat{G}}\{f_n^{p_1,\dots,p_{s-1},p_{s}}(t)\}_{n=0}^\infty,\ t\in[a,a+2l],\ p_1,\dots,p_s=1,\dots,k, 
\ee 
for every fixed $k,s\in\mbb{N}$.  

Sets (\ref{1.1}) -- (\ref{1.3}) give consequently the first generation, the second generation and so on up $s^{\mbox{th}}$ generation of new $L_2$-orthogonal systems. The counts of members of generations form the geometric sequence 
\bdis 
k,k^2,\dots,k^s. 
\edis  

\begin{remark}
If, for example, 
\bdis 
k=10^2,\ s=5\times 10^4, 
\edis  
i. e. for the $50\ 000^{\mbox{th}}$ generation of new $L_2$-orthogonal systems, we obtain 
\bdis 
10^{10^5}\ \mbox{multiple} 
\edis  
from one and only fixed $L_2$-orthogonal system. 
\end{remark} 

\subsection{} 

Let us remind the definition of the Legendre's polynomials 
\be \label{1.4} 
\{P_n(t)\}_{n=0}^\infty,\ t\in [-1,1] 
\ee  
by means of the generating function that is by the formula 
\be \label{1.5} 
\frac{1}{\sqrt{1-2ut-t^2}}=\sum_{n=0}^\infty P_n(t)u^n,\ u\in (-1,1), 
\ee  
where, of course, the collection (\ref{1.4}) represents the simple $L_2$-orthogonal system. 

Now, let us have a look on how the operator $\hat{G}$ acts on the system (\ref{1.4}). For given natural number $k$ the operator $\hat{G}$ produces as many as $k^3$ new species of $L_2$-orthogonal systems of the third generation as follows: 
\be \label{1.6} 
\begin{split}
& P_{n}^{p_1,p_2,p_3}(t)=P_n(u_{p_1}(u_{p_2}(u_{p_3}(t))))\times \prod_{r=0}^{p_1-1}|\tilde{Z}(v_{p_1}^r(u_{p_2}(u_{p_3}(t))))|\times \\ 
& \prod_{r=0}^{p_2-1}|\tilde{Z}(v_{p_2}^r(u_{p_3}(t)))|\times \prod_{r=0}^{p_3-1}|\tilde{Z}(v_{p_3}^r(t))|, \\ 
& p_1,p_2,p_3=1,\dots,k,\ t\in[-1,1],\ a=-1,\ l=1, 
\end{split}
\ee  
where 
\be \label{1.7} 
\begin{split}
& u_{p_i}(t)=\vp_1^{p_i}\left(\frac{\overset{p_i}{\wideparen{T+2}}-\overset{p_i}{T}}{2}(t+1)+\overset{p_i}{T}\right)-T-1,\ i=1,2,3, \\ 
& v_{p_1}^r(t)=\vp_1^{r}\left(\frac{\overset{p_i}{\wideparen{T+2}}-\overset{p_i}{T}}{2}(t+1)+\overset{p_i}{T}\right),\ r=0,1,\dots,p_i-1, \\ 
& t\in[-1,1] \ \Rightarrow\ u_{p_i}(t)\in [-1,1] \ \wedge \ v_{p_i}^r(t)\in [\overset{p_i-r}{T},\overset{p_i-r}{\wideparen{T+2}}]. 
\end{split}
\ee 

\subsection{} 

Now we give the following. 

\begin{mydef81} 
\begin{itemize}
	\item[(a)] Every member of every new $L_2$-orthogonal system 
	\bdis 
	\{ P_n^{p_1,p_2,p_3}(t)\}_{n=0}^{\infty},\ t\in[-1,1],\ p_1,p_2,p_3=1,\dots,k 
	\edis  
	contains the function 
	\bdis 
	\left|\zf\right|_{t=\tau} 
	\edis 
	for corresponding $\tau$ since (comp. \cite{3}, (9.1), (9.2)) 
	\be \label{1.8} 
	|\tilde{Z}(t)|=\sqrt{\frac{{\rm d}\vp_1(t)}{{\rm d}t}}=\{1+o(1)\}\frac{1}{\sqrt{\ln t}}|\zf|,\ t\to\infty. 
	\ee 
	\item[(b)] property (a) holds true due to the Theorem of this paper for every generation 
	\bdis 
	\{f_n^{p_1,\dots,p_s}(t)\}_{n=0}^\infty,\ t\in [a,a+2l],\ s\in\mbb{N}. 
	\edis 
\end{itemize}
\end{mydef81}  

\begin{remark}
The main aim of this paper is expressed by the Property 1. Namely, that there is a close binding between the theory of the Riemann's zeta-function on the critical line and the theory of $L_2$-orthogonal systems. Let us notice also that this paper finishes preparatory papers \cite{5} and \cite{6}. 
\end{remark} 

\section{Main result} 

\subsection{} 

We use the following notions: 
\begin{itemize}
	\item[(a)] Jacob's ladder $\vp_1(t)$, 
	\item[(b)] the function 
	\be \label{2.1} 
	\begin{split}
		& \tilde{Z}^2(t)=\frac{{\rm d}\vp_1(t)}{{\rm d}t}=\frac{1}{\omega(t)}\left|\zf\right|^2, \\ 
		& \omega(t)=\left\{1+\mcal{O}\left(\frac{\ln\ln t}{\ln t}\right)\right\}\ln t,\ t\to\infty, 
	\end{split}
	\ee  
	\item[(c)] direct iterations of the Jacob's ladder 
	\be \label{2.2} 
	\begin{split}
	& \vp_1^0(t)=t,\ \vp_1^1(t)=\vp_1(t),\ \vp_1^2(t)=\vp_1(\vp_1(t)), \dots, \vp_1^k(t)=\vp_1(\vp_1^{k-1}(t))
	\end{split}
	\ee 
	for every fixed $k\in\mbb{N}$, 
	\item[(d)] reverse iterations (by means of $\vp_1^{-1}(t)$) 
	\bdis 
	\begin{split}
	& [\onT,\overset{0}{\wideparen{T+U}}], [\ojT,\overset{1}{\wideparen{T+U}}],\dots, [\overset{k}{T},\overset{k}{\wideparen{T+U}}], \\ 
	& U=o\left(\frac{ T}{\ln T}\right),\ T\to\infty 
	\end{split}
	\edis 
	of the basic segment 
	\bdis 
	[T,T+U]=[\onT,\overset{0}{\wideparen{T+U}}], 
	\edis 
	that we have introduced into the theory of the Riemann's zeta-function, see \cite{1} -- \cite{5}. 
\end{itemize} 

\subsection{} 

Next we use the following analytic properties of Jacob's ladder $\vp_1(t)$: 
\begin{itemize}
	\item[(e)] $\vp_1(t)\in C^\infty([T_0,\infty])$ and it is strongly increasing function, see \cite{1}, 
	\item[(f)] $\vp_1^p(t)\in C^\infty([T_0,\infty]),\ p=1,\dots,k$ and it is again strongly increasing (this property follows easily from (e)) 
	\item[(g)] as a consequence of (f) we have next: every function 
	\bdis 
	\vp_1^p(t),\ t\in [A,B],\ A>T_0 
	\edis 
	is absolutely continuous and strongly increasing on every segment $[A,B]$ with $A>T_0$ (namely, the Lipschitz condition holds true by (f)), 
	\item[(h)] the composite function $F[f(t)]$, where $F$ is absolutely continuous and $f$ is absolutely continuous and monotonic is again absolutely continuous function. 
\end{itemize} 

\subsection{} 

Finally, we introduce the following functions (comp. (\ref{1.7})) together with some of their properties: 
\be \label{2.3} 
\begin{split}
& u_{p_i}(t)=\vp_1^{p_i}\left(\frac{\overset{p_i}{\wideparen{T+2l}}-\overset{p_i}{T}}{2l}(t-a)+\overset{p_i}{T}\right)-T+a, \\ 
& v_{p_i}^r(t)=\vp_1^{r}\left(\frac{\overset{p_i}{\wideparen{T+2l}}-\overset{p_i}{T}}{2l}(t-a)+\overset{p_i}{T}\right), \\ 
& t\in [a,a+2l],\ i=1,\dots,s,\ r=0,1,\dots,p_{i-1},\ p_i=1,\dots,k, 
\end{split}
\ee 
where 
\be \label{2.4} 
u_{p_i}(t)\in [a,a+2l],\ v_{p_i}^r(t)\in \left[ \overset{p_i-r}{T},\overset{p_i-r}{\wideparen{T+2l}}\right]
\ee  
and, with regard to the second inclusion, see \cite{5}, Property 2, the segments 
\bdis 
\left[ \overset{p_i-r}{T},\overset{p_i-r}{\wideparen{T+2l}}\right] 
\edis 
represent corresponding components of the disconnected set (see \cite{5}, (2.9)) 
\be \label{2.5} 
\Delta(T,k,l)=\bigcup_{r=0}^{k}\left[ \overset{r}{T},\overset{r}{\wideparen{T+2l}}\right] 
\ee  
and the following properties of the above mentioned set hold true\footnote{\cite{5}, (2.5), (2.6)}: 
\be \label{2.6} 
l=o\left(\frac{T}{\ln T}\right),\ T\to \infty \ \Rightarrow 
\ee 
\be \label{2.7} 
\left|\left[ \overset{r}{T},\overset{r}{\wideparen{T+2l}}\right] \right|=\overset{r}{\wideparen{T+2l}}-\overset{r}{T}=o\left(\frac{T}{\ln T}\right), 
\ee 
\be \label{2.8} 
\left|\left[ \overset{r-1}{\wideparen{T+2l}},\overset{r}{T}\right] \right|=\overset{r}{T}-\overset{r-1}{\wideparen{T+2l}}\sim (1-c)\frac{T}{\ln T}, 
\ee 
\be \label{2.9} 
\left[ \overset{0}{T},\overset{0}{\wideparen{T+2l}}\right]\prec \left[ \overset{1}{T},\overset{1}{\wideparen{T+2l}}\right]\prec \dots \prec \left[ \overset{k}{T},\overset{k}{\wideparen{T+2l}}\right], 
\ee  
where $c$ is the Euler's constant and the property (\ref{2.9}) follows from (\ref{2.8}). 

\begin{remark}
Asymptotic behavior of the disconnected set (\ref{2.5}) is as follows: if $T\to\infty$, then the components of this set recede unboundedly each from other and all together are receding to infinity. Hence the set (\ref{2.5}) behaves at $T\to\infty$ as one-dimensional Friedmann-Hubble expanding universe. 
\end{remark} 

\subsection{} 

The following Theorem is the main result of this paper.  

\begin{mydef11}
There is such a generating vector-operator $\hat{G}$ that for every fixed $L_2$-orthogonal system 
\be \label{2.10} 
\{f_n(t)\}_{n=0}^\infty,\ t\in [a,a+2l],\ a\in \mbb{R}, l\in\mbb{R}^+ 
\ee 
and for every fixed $k\in\mbb{N}$ the operator $\hat{G}$ associates following orthogonal systems: 
\begin{itemize}
	\item[(a)] the first generation of the following new species of $L_2$-orthogonal systems 
	\be \label{2.11} 
	\{f_n^{p_1}(t)\}_{n=0}^\infty,\ t\in [a,a+2l],\ p_1=1,\dots,k
	\ee 
	i. e. 
	\be \label{2.12} 
	\hat{G}[\{f_n(t)\}]=(\{f_n^1(t)\},\{f_n^2(t)\},\dots,\{f_n^{k}(t)\});\ \{f_n(t)\}=\{f_n(t)\}_{n=0}^\infty, 
	\ee 
	where\footnote{Comp. (\ref{2.3}).} 
	\be \label{2.13} 
	f_n^{p_1}(t)=f_n(u_{p_1}(t))\prod_{r=0}^{p_1-1}\left|\tilde{Z}(v_{p_1}^r(t))\right| 
	\ee  
	and every of the functions $u_{p_1}(t)$ defines an automorphism on $[a,a+2l]$, 
	\item[(b)] the second generation of the following $L_2$-orthogonal systems 
	\be \label{2.14} 
	\{f_n^{p_1,p_2}(t)\}_{n=0}^\infty,\ t\in [a,a+2l],\ p_1,p_2=1,\dots,k, 
	\ee 
	i. e. 
	\be \label{2.15} 
	\hat{G}[\{f_n^{p_1}(t)\}]=(\{f_n^{p_1,1}(t)\},\{f_n^{p_1,2}(t)\},\dots,\{f_n^{p_1,k}(t)\}),  
	\ee 
	where 
	\be \label{2.16} 
	f_n^{p_1,p_2}(t)=f_n(u_{p_1}(u_{p_2}(t)))\prod_{r=0}^{p_1-1}\left|\tilde{Z}(v_{p_1}^r(u_{p_2}(t)))\right|\prod_{r=0}^{p_2-1}\left| \tilde{Z}(v_{p_2}^r(t))\right| 
	\ee  
	and each of the functions $u_{p_1}(u_{p_2}(t))$ defines an automorphism on $[a,a+2l]$, 
	\item[(c)] and so on up to the $s^{\mbox{th}}$ geneation of the $L_2$-orthogonal systems 
	\be \label{2.17} 
	\{ f_n^{p_1,p_2,\dots,p_s}(t)\}_{n=0}^{\infty},\ t\in [a,a+2l],\ p_1,\dots,p_s=1,\dots,k, 
	\ee  
	i. e. 
	\be \label{2.18} 
	\hat{G}[\{f_n^{p_1,\dots,p_{s-1}}(t)\}]=(\{f_n^{p_1,\dots,p_{s-1},1}(t)\},\{f_n^{p_1,\dots,p_{s-1},2}(t)\},\dots,\{f_n^{p_1,\dots,p_{s-1},k}(t)\}), 
	\ee  
	where 
	\be \label{2.19} 
	\begin{split}
	& f_n^{p_1,p_2,\dots,p_s}(t)=f_n(u_{p_1}(u_{p_2}(\dots (u_{p_s}(t))\dots)))\times \\ 
	& \prod_{r=0}^{p_1-1}\left|\tilde{Z}(v_{p_1}^r(u_{p_2}(u_{p_3}(\dots (u_{p_s}(t))\dots ))))\right|\times  \\ 
	& \prod_{r=0}^{p_2-1}\left|\tilde{Z}(v_{p_2}^r(u_{p_3}(u_{p_3}(\dots (u_{p_s}(t))\dots ))))\right|\times \\ 
	& \vdots \\ 
	& \prod_{r=0}^{p_{s-1}-1}\left|\tilde{Z}(v_{p_s-1}^r(u_{p_s}(t)))\right| \times \\ 
	& \prod_{r=0}^{p_s-1}\left|\tilde{Z}(v_{p_s}(t))\right| 
	\end{split}
	\ee 
	and each of the functions 
	$$u_{p_1}(u_{p_2}(\dots (u_{p_s}(t))\dots ))$$ 
	defines an automorphism on $[a,a+2l]$, 
	\item[(d)] for every fixed $k,s\in\mbb{N}$ the $L_2$-orthonormal system 
	\be \label{2.20} 
	\left(\prod_{i=0}^s\sqrt{\frac{2l}{\overset{i}{\wideparen{T+2l}}-\overset{i}{T}}}\right)f_n^{p_1,\dots,p_s}(t),\ t\in[a,a+2l],\ p_1,\dots,p_s=1,\dots,k
	\ee 
	is corresponding with (\ref{2.19}), 
	\item[(e)] finally, all these formulas are true for all sufficiently big $T>0$, that is we have the continuum set of possibilities how to construct new classes of $L_2$-orthogonal systems (\ref{2.17}). 
\end{itemize}
\end{mydef11}

\subsection{} 

Let us denote the set of all $L_2$-orthogonal systems not containing the functions 
\be \label{2.21} 
\vp_1(t),\ \tilde{Z}^2(t)=\frac{{\rm d}\vp_1(t)}{{\rm d}t}
\ee  
as 
\bdis 
L_2^0([a,a+2l]). 
\edis 
Now it is true by Theorem 1 that for every fixed 
\bdis 
\begin{split}
& \{f_n(t)\}_{n=0}^\infty = \{f_n(t)\}\in L_2^0([a,a+2l]), \\ 
& \{f_n(t)\}\xrightarrow{\hat{G}}\{ f_n^{p_1}(t)\}\xrightarrow{\hat{G}}\dots \xrightarrow{\hat{G}} \{ f_n^{p_1,\dots,p_s}(t)\}. 
\end{split}
\edis 

We instantly get the following. 

\begin{mydef41}
\be \label{2.22} 
L_2^0([a,a+2l])\xrightarrow{\hat{G}}L_2^1([a,a+2l])\xrightarrow{\hat{G}}\dots \xrightarrow{\hat{G}} L_2^s([a,a+2l]), 
\ee  
where 
\bdis 
L_2^i([a,a+2l])
\edis 
stands for the $i^{\mbox{th}}$ generation of the image of $L_2^0([a,a+2l])$. Consequently, the union 
\be \label{2.23} 
\bigcup_{i=1}^s L_2^i([a,a+2l])
\ee  
represents a kind of $\zeta$-extension of the $s^{\mbox{th}}$ order of $L_2^0([a,a+2l])$ for every fixed segment $[a,a+2l]$ and every fixed $k,s\in\mbb{N}$. Finally, the union over all the segments 
\be \label{2.24} 
\bigcup_{a\in \mbb{R}, l>0}\left\{ \bigcup_{i=1}^s L_2^i([a,a+2l])\right\}
\ee 
represents the complete $\zeta$-extension of the $s^{\mbox{th}}$ order of the set 
\be \label{2.25} 
\bigcup_{a\in \mbb{R}, l>0} L_2^0([a,a+2l]). 
\ee 
\end{mydef41} 

\subsection{} 

\begin{remark}
We may select another function 
\bdis 
\psi(t)\not=\vp_1(t),\ \psi(t)\nsim \vp_1(t),\ t\to\infty 
\edis 
instead of the Jacob's ladder and use the pair 
\bdis 
\psi(t),\ \frac{{\rm d}\psi(t)}{{\rm d}t}
\edis 
to extend $L_2^0([a,a+2l])$. This way however will not attain $\zeta$-extensions (\ref{2.13}), (\ref{2.16}), (\ref{2.19}), (\ref{2.23}) and (\ref{2.24}) and therefore will be irrelevant. 
\end{remark}  

\section{Jacob's ladders} 

\subsection{} 

Let us remind that the Jacob's ladder 
\bdis 
\vp_1(t)=\frac{1}{2}\vp(t) 
\edis  
was introduced in \cite{1}, see also \cite{3}, where the function $\vp(t)$ is an arbitrary continuous solution of the nonlinear integral equation\footnote{Also introduced in \cite{1}.} 
\be \label{3.1} 
\int_0^{\mu[x(T)]}Z^2(t)e^{-\frac{2}{x(T)}t}{\rm d}t=\int_0^TZ^2(t){\rm d}t, 
\ee  
where 
\be \label{3.2} 
\begin{split}
& Z(t)=e^{i\vartheta(t)}\zf, \\ 
& \vartheta(t)=-\frac t2\ln\pi+\im\left\{\ln\Gamma\left(\frac 14+i\frac t2\right)\right\} 
\end{split} 
\ee 
and the class of functions $\{\mu\}$ is specified as 
\bdis 
\mu\in C^\infty([y_0,+\infty))
\edis 
being monotonically increasing, unbounded from above and obeying the inequality 
\be \label{3.3} 
\mu(y)\geq 7y\ln y. 
\ee 
Every admissible function $\mu(y)$ generates a solution 
\bdis 
y=\vp(T;\mu)=\vp(T). 
\edis   

\begin{remark}
The function $\vp_1(T)$ is called Jacob's ladder as an analogue of the Jacob's dream in Chumash, Bereishis, 28:12. 
\end{remark} 

\subsection{} 

Let us remind that the Hardy-Littlewood integral (1918) 
\be \label{3.4} 
\int_0^T\left|\zf\right|^2{\rm d}t 
\ee  
can be expressed as follows: 
\be \label{3.5} 
\int_0^T\left|\zf\right|^2{\rm d}t=T\ln T+(2c-1-\ln 2\pi)T+R(T), 
\ee  
with, for example, Ingham's error term 
\be \label{3.6} 
R(T)=\mcal{O}(T^{1/2}\ln T)=\mcal{O}(T^{1/2+\delta}),\ \delta>0,\ T\to\infty 
\ee 
for arbitrary small $\delta$. 

Next, it is true by Good's $\Omega$-theorem (1977), that 
\be \label{3.7} 
R(T)=\Omega(T^{1/4}),\ T\to\infty. 
\ee  

\begin{remark}
Let 
\be \label{3.8} 
R_a(T)=\mcal{O}(T^{1/4+a}),\ a\in \left[ \delta,\frac 14+\delta\right],\ T\to\infty. 
\ee 
Then, by (\ref{3.7}), one obtain for every valid estimate of the type (\ref{3.8}) that 
\be\label{3.9} 
\limsup_{T\to\infty}|R_a(T)|=+\infty. 
\ee 
In other words, every expression of the type (\ref{3.5}) and (\ref{3.9}) possesses an unbounded error at infinity. 
\end{remark} 

\subsection{} 

Under the circumstances (\ref{3.5}) and (\ref{3.9}) we have shown in our paper \cite{1} that the Hardy-Littlewood integral (\ref{3.4}) has an infinite set of almost exact representations expressed by the following formula. 

\begin{mydef91}
\be \label{3.10} 
\begin{split}
& \int_0^T\left|\zf\right|^2{\rm d}t= \\ 
& \vp_1(T)\ln\{\vp_1(T)\}+(c-\ln 2\pi)\vp_1(T)+c_0+\mcal{O}\left(\frac{\ln T}{T}\right),\ T\to\infty 
\end{split}
\ee 
($c$ is the Euler's constant and $c_0$ is the constant from the Titchmarsh-Kober-Atkinson formula) with the error term vanishing at infinity: 
\be \label{3.11} 
\tilde{R}(T)=\mcal{O}\left(\frac{\ln T}{T}\right)\xrightarrow{T\to\infty}0. 
\ee   
\end{mydef91} 

\begin{remark}
The comparison of (\ref{3.9}) and (\ref{3.11}) completely characterizes the level of exactness of our representation (\ref{3.10}) of the Hardy-Littlewood integral (\ref{3.4}). 
\end{remark} 

\subsection{} 

Now, let us remind other formulae demonstrating the power of Jacob's ladder $\vp_1(t)$. 

First we have obtained the following\footnote{See \cite{3}, (8.3).}  

\begin{mydef92}
\be \label{3.12} 
\begin{split}
& \int_T^{T+U}\left|\zfvp\right|^4\left|\zf\right|^2{\rm d}t\sim \frac{1}{2\pi^2}U\ln^5T, \\ 
& U=T^{7/8+2\delta},\ T\to\infty. 
\end{split}
\ee 
\end{mydef92} 

\begin{remark}
The formula (\ref{3.12}) is the first asymptotic formula of the sixth order in $|\zeta|$ on the critical line in the theory of the Riemann's zeta-function. 
\end{remark}  

Next, let 
\bdis 
S(t)=\frac{1}{\pi}\arg\left\{\zf\right\},\ S_1(T)=\int_0^TS(t){\rm d}t, 
\edis  
where the function $\arg$ is defined in the usual way. We have obtained the following two formulae concerning $S(t)$ function:\footnote{See \cite{2}, (5.4), (5.5).} 

\begin{mydef93}
\be \label{3.13} 
\begin{split}
& \int_T^{T+U}\left[\arg\left\{\zfvp\right\}\right]^{2k}\left|\zf\right|^2{\rm d}t\sim \\ 
& \frac{1}{2^k}U\ln T(\ln\ln T)^k,\ U\in \left[ T^{1/3+\delta},\frac{T}{\ln T}\right],\ T\to\infty, 
\end{split}
\ee 
for every fixed $k\in\mbb{N}$. 
\end{mydef93}  

\begin{mydef94}
\be \label{3.14} 
\int_T^{T+U}\{S_1[\vp_1(t)]\}^{2k}\left|\zf\right|^2{\rm d}t\sim a_kU\ln T,\ T\to\infty. 
\ee 
\end{mydef94} 

\begin{remark}
New kind of the classical A. Selberg's formulae (1946) are expressed by means of our results (\ref{3.13}) and (\ref{3.14}). 
\end{remark} 

\begin{remark}
Let us notice explicitly, that formulae (\ref{3.12}) -- (\ref{3.14}) are $\zeta$-correlation formulae on the critical line. For example, (\ref{3.13}) at $k=1$ describes interaction between values of the functions 
\bdis 
\left[\arg\left\{\zfvp\right\}\right]^{k},\ \left|\zf\right|^2. 
\edis 
\end{remark} 

\section{Proof of Theorem 1} 

\subsection{} 

In \cite{5}, (7.1) (7.2), we have shown the following results concerning direct and reverse iterations: 

\begin{mydef51}
If 
\be \label{4.1} 
U=o\left(\frac{T}{\ln T}\right),\ T\to\infty, 
\ee  
then for every function 
\bdis 
g(t)\in L([T,T+U]) 
\edis  
the following holds true: 
\be \label{4.2} 
\int_T^{T+U}g(t){\rm d}t=\int_{\overset{p}{T}}^{\overset{p}{\wideparen{T+U}}}g[\vp_1^p(\tau)]\prod_{r=0}^{p-1}\tilde{Z}^2[\vp_1^r(\tau)]{\rm d}\tau,\ 
p=1,\dots,k
\ee  
for every fixed $k\in\mbb{N}$. 
\end{mydef51} 

\begin{remark}
Let us notice also the subsection 2.2, (g). 
\end{remark} 

\subsection{Proof of Theorem 1} 

\subsubsection{} 

Let 
\be \label{4.3} 
\{f_n(t)\}_{n=0}^\infty \subset L_2([a,a+2l])
\ee  
is arbitrary fixed system of orthogonal functions. Then $l$ is also fixed positive number and condition (\ref{4.1}) is fulfilled for $U=2l$ for all sufficiently big and positive $T$. Now we have ($m\not=n$) 
\be \label{4.4} 
0=\int_a^{a+2l}f_m(t)f_n(t){\rm d}t=\int_T^{T+2l}f_m(\tau-T+a)f_n(\tau-T+a){\rm d}\tau =  
\ee 
next we obtain by Lemma 1 for any sufficiently big $T$ 
\be \label{4.5} 
 = \int_{\overset{p}{T}}^{\overset{p}{\wideparen{T+2l}}}f_m[\vp_1^p(\rho)-T+a]f_n[\vp_1^p(\rho)-T+a]\prod_{r=0}^{p-1}\tilde{Z}^2[\vp_1^r(\rho)]{\rm d}\rho=
\ee 
and next, by simple sunstitution 
\bdis 
\rho=\rho(t)=\frac{\overset{p}{\wideparen{T+2l}}-\overset{p}{T}}{2l}(t-a)+\overset{p}{T},\ \rho\in [\overset{p}{T},\overset{p}{\wideparen{T+2l}}],\ t\in[a,a+2l], 
\edis  
where $\rho(t)$ is absolutely continuous and increasing, we obtain 
\be \label{4.6} 
\begin{split}
& = \frac{\overset{p}{\wideparen{T+2l}}-\overset{p}{T}}{2l}\int_a^{a+2l}f_m\left[\vp_1^p\left(\frac{\overset{p}{\wideparen{T+2l}}-\overset{p}{T}}{2l}(t-a)+\overset{p}{T}\right)-T+a\right]\times \\ 
& f_n\left[\vp_1^p\left(\frac{\overset{p}{\wideparen{T+2l}}-\overset{p}{T}}{2l}(t-a)+\overset{p}{T}\right)-T+a\right]\times \\ 
& \prod_{r=0}^{p-1}\tilde{Z}^2\left[\vp_1^r\left(\frac{\overset{p}{\wideparen{T+2l}}-\overset{p}{T}}{2l}(t-a)+\overset{p}{T}\right)\right]{\rm d}t=
\end{split}
\ee 
and, in the next step of the first cycle, we finish with 
\be \label{4.7} 
\begin{split}
& = \frac{\overset{p}{\wideparen{T+2l}}-\overset{p}{T}}{2l}\int_a^{a+2l}f_m^p(t)f_n^p(t){\rm d}t \ \Rightarrow \int_a^{a+2l}f_m^p(t)f_n^p(t){\rm d}t=0, 
\end{split}
\ee 
where 
\be \label{4.8} 
\begin{split}
& f_n^p(t)=f_n\left[\vp_1^p\left(\frac{\overset{p}{\wideparen{T+2l}}-\overset{p}{T}}{2l}(t-a)+\overset{p}{T}\right)-T+a\right]\times \\ 
& \prod_{r=0}^{p-1}\left|\tilde{Z}\left[\vp_1^r\left(\frac{\overset{p}{\wideparen{T+2l}}-\overset{p}{T}}{2l}(t-a)+\overset{p}{T}\right)\right]\right|,\ t\in [a,a+2l],\ p=1,\dots,k. 
\end{split}
\ee 

\subsubsection{} Now we give the following 

\begin{mydef21} 
The symbol $\hat{G}$ stands for vector operator defined on the set all $L_2$-orthogonal systems 
\bdis 
\{f_n(t)\}_{n=0}^\infty,\ t\in [a,a+2l],\ a\in\mbb{R},\ l>0
\edis 
defined by three integral transformatrions (\ref{4.4}) -- (\ref{4.6}). $\hat{G}$ maps an $L_2$-orthogonal system into $k$-tuple of new orthogonal systems 
\bdis 
\begin{split}
& \hat{G}[\{f_n(t)\}]=(\{f_n^1(t)\},\{f_n^2(t)\},\dots,\{f_n^k(t)\})=\{f_n^p(t)\}_{n=0}^\infty, \\ 
& p=1,\dots,k,\ t\in [a,a+2l];\ \{f_n(t)\}_{n=0}^\infty = \{f_n(t)\},\dots 
\end{split} 
\edis 
for every fixed $k\in\mbb{N}$. 
\end{mydef21} 

\subsubsection{} 

Let us notice that the transformation\footnote{Comp. (\ref{2.3}), (\ref{4.8}).} 
\be \label{4.9} 
u_p(t)=\vp_1^p\left(\frac{\overset{p}{\wideparen{T+2l}}-\overset{p}{T}}{2l}(t-a)+\overset{p}{T}\right)-T+a,\ t\in [a,a+2l] 
\ee  
has the following properties: 
\begin{itemize}
	\item[(a)] By \cite{5}, subsection 6.1 
	\bdis 
	\begin{split}
	& u_p(a)=\vp_1^p(\overset{p}{T})-T+a=T-T+a=a;\ \overset{p}{T}=\vp_1^{-p}(T), \\ 
	& u_p(a+2l)=\vp_1(\overset{p}{\wideparen{T+2l}})-T+a=a+2l. 
	\end{split}
	\edis 
	\item[(b)] Since the continuous function $\vp_1^p(t)$ is increasing and 
	\bdis 
	\rho=\frac{\overset{p}{\wideparen{T+2l}}-\overset{p}{T}}{2l}(t-a)+\overset{p}{T},\ t\in[a,a+2l] 
	\edis 
	it is evident that the composite function 
	\bdis 
	u_p(t),\ t\in[a,a+2l] 
	\edis 
	is also increasing and therefore 
	\bdis 
	t\in[a,a+2l] \ \Rightarrow \ u_p(t)\in [a,a+2l]. 
	\edis 
\end{itemize} 

\begin{remark}
We have as a consequence of (a) and (b) that new automorphism on $[a,a+2l]$ is defined by the one-to-one correspondence (\ref{4.9}) for every fixed sufficiently big positive $T$. Of course, every function $u_{p_i}(t)$ defines an automorphism on $[a,a+2l]$ too. 
\end{remark} 

\subsubsection{} 

By making use of the operator $\hat{G}$ on the system 
\bdis 
\{f_n^{p_1}(t)\}_{n=0}^\infty,\ t\in [a,a+2l],\ p_1=1,\dots,k 
\edis 
(the second cycle) we obtain 
\be \label{4.10} 
\hat{G}[\{f_n^{p_1}(t)\}]=\{f_n^{p_1,p_2}(t)\}_{n=0}^\infty,\ p_1,p_2=1,\dots,k, 
\ee  
where 
\be \label{4.11} 
\begin{split}
& f_n^{p_1,p_2}(t) = \\ 
& f_n\left[\vp_1^{p_1}\left(\frac{\overset{p_1}{\wideparen{T+2l}}-\overset{p_1}{T}}{2l}\left(\vp_1^{p_2}\left(\frac{\overset{p_2}{\wideparen{T+2l}}-\overset{p_2}{T}}{2l}(t-a)+\overset{p_2}{T}\right)-T\right)+\overset{p_1}{T}\right)-T+a\right]\times \\ 
& \prod_{r=0}^{p_1-1}\left|\tilde{Z}\left[\vp_1^r\left(\left(\frac{\overset{p_1}{\wideparen{T+2l}}-\overset{p_1}{T}}{2l}\left(\vp_1^{p_2}\left(\frac{\overset{p_2}{\wideparen{T+2l}}-\overset{p_2}{T}}{2l}(t-a)+\overset{p_2}{T}\right)-T\right)+\overset{p_1}{T}\right)\right)\right]\right|\times \\ 
& \prod_{r=0}^{p_2-1}\left|\tilde{Z}\left[\vp_1^r\left(\frac{\overset{p_2}{\wideparen{T+2l}}-\overset{p_2}{T}}{2l}(t-a)+\overset{p_2}{T}\right)\right]\right|. 
\end{split}
\ee 
It is clear that there is need for simplification of our formulae (\ref{4.8}) and (\ref{4.11}). For this purpose we use functions of (\ref{2.3})\footnote{See also Remark 11.} that provides us with the results 
\be \label{4.12} 
\begin{split}
& f_n^{p_1}(t)=f_n(u_{p_1}(t))\prod_{r=0}^{p_1-1}|\tilde{Z}(v_{p_1}^r(t))|, \\ 
& f_n^{p_1,p_2}(t)=f_n(u_{p_1}(u_{p_2}(t)))\prod_{r=0}^{p_1-1}|v_{p_1}^r(u_{p_2}(t))|\prod_{r=0}^{p_2-1}|\tilde{Z}(v_{p_2}^r(t))|, 
\end{split}
\ee 
where, for example, 
\bdis 
\{u_{p_1}(t_1),\ t_1=u_{p_2}(t_2),\ t_1,t_2\in [a,a+2l]\}\ \Rightarrow \ u_{p_1}(u_{p_2}(t_2)),\ t_2=t, 
\edis  
i. e. we have formulae (\ref{2.13}) and (\ref{2.16}). 

\subsubsection{} 

Next, in the $(s-1)^{\mbox{th}}$-cycle we get 
\be \label{4.13}  
\begin{split}
& f_n^{p_1,p_2,\dots,p_{s-1}}(t)=f_n(u_{p_1}(u_{p_2}(\dots (u_{p_{s-1}}(t))\dots)))\times \\ 
& \prod_{r=0}^{p_1-1}\left|\tilde{Z}(v_{p_1}^r(u_{p_2}(u_{p_3}(\dots (u_{p_{s-1}}(t))\dots ))))\right|\times  \\ 
& \prod_{r=0}^{p_2-1}\left|\tilde{Z}(v_{p_2}^r(u_{p_3}(u_{p_4}(\dots (u_{p_{s-1}}(t))\dots ))))\right|\times \\ 
& \vdots \\ 
& \prod_{r=0}^{p_{s-2}-1}\left|\tilde{Z}(v_{p_s-2}^r(u_{p_{s-1}}(t)))\right| \times \\ 
& \prod_{r=0}^{p_{s-1}-1}\left|\tilde{Z}(v_{p_{s-1}}(t))\right|,\ s>2, 
\end{split}
\ee 
then, if we use the operator $\hat{G}$ on (\ref{4.13}) to obtain the $s^{\mbox{th}}$ cycle, we get the set of $k^s$ formulas (\ref{2.19}). That means the formula (\ref{2.15}) holds true for every $s\in\mbb{N}$. 

I would like to thank Michal Demetrian for his moral support of my study of Jacob's ladders.

\end{document}